\documentclass[journal]{IAENGtran}
\ifCLASSINFOpdf
   \usepackage[pdftex]{graphicx}
   \DeclareGraphicsExtensions{.pdf,.jpeg,.png}
\else
   \usepackage[dvips]{graphicx}
   \DeclareGraphicsExtensions{.eps}
\fi
\usepackage{amssymb}
\usepackage{mathtools}

\begin{document}
%
\title{A study of the inverse scattering problem for the reconstruction of the shape and/or impedance of an obstacle}
%
%
%

\author{Sarika~Karanth\IAENGmembership{},~Shobha~M~Erappa~\IAENGmembership{}
\thanks{Manuscript received September 22, 2022.}
\thanks{Sarika Karanth holds a Master’s in Applied Mathematics and Computing from the Department of Mathematics, Manipal Institute of Technology, Manipal Academy of Higher Education, India (e-mail: sarika.karanth@gmail.com). }
\thanks{Shobha M Erappa is an Assistant Professor – Selection Grade of the Department of Mathematics, Manipal Institute of Technology, Manipal Academy of Higher Education, India (e-mail: shobha.me@manipal.edu).}}

\maketitle

\pagestyle{empty}
\thispagestyle{empty}

\begin{abstract}
Three papers describing different methods to solve the inverse scattering problem of the reconstruction of the shape and/or impedance of an obstacle have been chosen for analysis. This literature review consists of an evaluation of these methods in which comparison of the assumptions, conditions, advantages, disadvantages, accuracy, and applicability has been undertaken. 
\end{abstract}

\begin{IAENGkeywords}
Helmholtz-equation, Impedance-reconstruction, Inverse-obstacle-scattering, Shape-reconstruction.
\end{IAENGkeywords}

%
\IAENGpeerreviewmaketitle

\section{Introduction}
%
%
%
%
\IAENGPARstart{I}{nverse} problems are those that use the examination of a response of an obstacle to a external signal to determine an unspecified property of that obstacle \cite{AGR}. In an inverse scattering problem, an unknown property of the object is found by directing an external acoustic or electromagnetic field onto the object and observing the scattered field \cite{AGR}. This problem has found many applications due to its non-invasiveness, some of which include optics, geophysics, medical imaging, non-destructive evaluation and seismic exploration \cite{EEH}.\par
In this review, we have selected three papers, namely Three-dimensional time harmonic electromagnetic inverse scattering: the reconstruction of the shape and the impedance of an obstacle \cite{CMA}, Numerical solution of an inverse obstacle scattering problem with near-field data \cite{JCP}, and Reconstruction of shapes and impedance functions using few far-field measurements \cite{RSI}. These papers describe different methods to solve the inverse scattering problem of reconstructing the shape and/or impedance function of an obstacle \cite{IP}. An evaluation of these methods has been conducted in this paper. 


\section{Three-dimensional time harmonic electromagnetic inverse scattering: the reconstuction of the shape and the impedance of an obstacle \cite{CMA}}
This paper illustrates a numerical method for the reconstruction of the boundary structure and impedance function of an obstacle from electromagnetic time harmonic scattering data.\par 
Let  $\mathbb{R}^3$ 
be the three-dimensional real Euclidean space and $x = (x_1,x_2,x_3)^\mathsf{T} \in  \mathbb{R}^3$ be a non-specific vector. Let $D$ $\subset$ $\mathbb{R}^3$ be a simply connected, bounded region with a smooth boundary $\partial D$ and it is assumed that $D$ contains the origin. The boundary of the obstacle has an electrical impedance given by $\chi (x)$, where $x \in \partial D$. $D$ is assumed to be in the shape of a star and symmetric with respect to the $x_3$-axis.\par
Let the equation for the total electric field be
\begin{equation} \label{eu_eqn}
{E(x) = E^i (x) + E^s (x)} 
\end{equation}
where $E^i (x) = \omega e^{ik(x,\alpha)}$ is the element which depends on the space variables $x$ of the electric field associated to a linearly-polarized incident plane wave, $E^s (x)$
is the element which depends on the space variables $x$ of the field scattered by the object when the incoming wave $E^i (x)$ is incident on it,  $\omega \in \mathbb{R}^3$ is the polarization vector, $\alpha \in \mathbb{R}^3$ is the direction of propagation of the plane wave, $\|\alpha\|$ = 1, $\omega$, $\alpha$ are provided and $k>0$ denotes the wavenumber.\par
The scattered field is given by, 
\begin{equation} \label{eu_eqn}
{E^s (x) = \frac{e^{ik\|x\|}}{\|x\|}E_0(\hat{x}, k, \alpha, \omega) + O\left(\frac{1}{\|x\|^2}\right)}
\end{equation}
${\|x\|}\rightarrow\infty$, where $E_0$ denotes the electric far field pattern which is created by the interaction of the object $D$ with $E^i (x)$, $\hat{x}=\Bigl(\frac{x}{\|x\|}\Bigl)$, $x\neq 0$.\par

Using the known information of the electrical far field patterns formed by many incident waves, the reconstruction of the boundary  $\partial D$ and the boundary electrical impedance $\chi (x)$ of $D$  is carried out. This method of reconstruction is derived from the Herglotz function method \cite{JSSC}.

\subsection{Governing equation}
The governing equation is the Helmholtz equation given by,
\begin{equation} \label{eu_eqn}
{\Delta E^s(x)+k^2E^s(x)=0, \   \mathbb{R}^3\backslash\bar{D}}
\end{equation}
with the divergence free condition,
\begin{equation} \label{eu_eqn}
{\text{div} E^s(x)=0,\ 
\mathbb{R}^3\backslash\bar{D}}
\end{equation}
where $\Delta$ is the vector Laplace operator.

\subsection{Assumptions and Conditions}
\begin{itemize}
  \item In this paper, the shape and impedance function are found.
  \item $\mathbb{R}^3$ is filled with a homogenous, isotropic medium.
  \item $D\subset \mathbb{R}^3$ is a simply-connected, bounded region with a smooth boundary $\partial D$ and it is assumed that the origin is contained within $D$.
  \item It is presumed that the incoming waves are plane waves and that the fields are electromagnetic in nature and time harmonic.
  \item Let ($r, \theta, \psi$) be the spherical coordinates of $\mathbb{R}^3$. The obstacle is presumed to be star shaped and symmetric with respect to the $x_3$ - axis, that is, 
  \begin{equation} \label{eu_eqn}
{\partial D = \{(r, \theta, \psi)\big|r = f(\theta),\  0\leq\theta\leq\pi\}}
\end{equation}
and the impedance is given by,
\begin{equation} \label{eu_eqn}
{\chi (\theta) = \chi ^{Re}(\theta)+i\chi ^{Im}(\theta),\  0\leq\theta\leq\pi}.
\end{equation}
The smooth nature of $\partial D$ provides an implication that $f$ is smooth and the assumption of symmetry of $\partial D$ implies,
\begin{equation} \label{eu_eqn}
{\frac{df}{d\theta}(0) = \frac{df}{d\theta}(\pi) = 0}.
\end{equation}
  \item $(\omega, \alpha) = 0$, that is, div $E^i (x) = 0$
  \item The governing equations are equipped with two boundary conditions which are as follows:
  \begin{itemize}
      \item[$*$] The Silver-Müller radiation condition is strictly related to the character of the electromagnetic fields and is at infinity. It is given by,
      \begin{equation} \label{eu_eqn}
{\text{curl} E^s (x) \times \hat{x} - ikE^s (x) = O\left(\frac{1}{\|x\|}\right),{\|x\|}\rightarrow\infty}
\end{equation}
where $\hat{x}=\Bigl(\frac{x}{\|x\|}\Bigl)$, $x\neq 0$.

      \item[$*$] \begin{equation} \label{eu_eqn}
{\hat{v}(x) \times \text{curl} E(x) + ik\chi (x) \hat{v}(x) \times (\hat{v}(x) \times  E(x)) = 0}
\end{equation}
      where $x \in \partial D$, $\hat{v}(x)$ is the exterior unit normal to $\partial D$ in the point $x \in \partial D$ and $\chi (x)$, $x \in \partial D$ is the boundary electrical impedance of the obstacle. This condition depends on the obstacle's electrical nature.
  \end{itemize}
  \item It is considered that if  $L$ is the characteristic length of the obstacle, then $kL\approx 1$, that is, the algorithm is effective in the resonance region.
\end{itemize}
\section{Numerical solution of an inverse obstacle scattering problem with near-field data \cite{JCP}}
In this paper, a generic, time-harmonic wave is incident to a sound-soft object, resulting in a scattered field. Here, the obstacle surface is reconstructed from the field measured on a circle surrounding the object by using near field data \cite{IPI}-\cite{CCP}. \par
The object is considered to be a small, smooth, perturbation of a disk in $\mathbb{R}^2$. In the polar coordinate, the obstacle is described by a domain, 
\begin{equation} \label{eu_eqn}
{\Omega = \{ (r, \theta): 0<r<a+f(\theta), \ \theta \in [0, 2\pi] \}}
\end{equation}
where $a>0$ denotes the radius of the unperturbed disk and $f(\theta)$ is the obstacle surface function and a $2\pi$ - periodic function with an infinity norm small in comparison to the wavelength $\lambda$ of the incoming wave. Hence, it is given by,
\begin{equation} \label{eu_eqn}
{f(\theta) = \varepsilon g (\theta)}
\end{equation}
where $g$ is the obstacle profile function such that $\|g\|_{\infty} = O(\lambda)$ and $\varepsilon > 0$ is a small perturbation parameter. $\mathbb{R}^2 \backslash \bar{\Omega}$
is assumed to be filled with a homogenous medium characterized by a constant wavenumber $k = \frac{2\pi}{\lambda}$.\par
The circle on which the total field $u$  is measured is given by,
\begin{equation} \label{eu_eqn}
{\Gamma = \{ (b, \theta): r = b>a+\|f\|_{\infty}, \ \theta \in [0, 2\pi] \}}.
\end{equation}\par
An incident wave $u^{inc} (r,\theta)$ on the obstacle produces a scattered wave given by $u^{sca}$. Due to the obstacle being sound-soft, the total field on the surface of the obstacle becomes null, that is,  
\begin{equation} \label{eu_eqn}
{u = 0}.
\end{equation}
Given the incident field $u^{inc}$, the inverse problem \cite{CK} is the reconstruction of the surface function $f(\theta)$ from the total field $u$ measured on $\Gamma$.

\subsection{Governing equation}
The incident field is denoted by $u^{inc}$ and the scattered field, by $u^{sca}$. \par
The governing equation is the Helmholtz equation. The scattered field $u^{sca}$ satisfies this equation and hence, the following equation is inferred,
\begin{equation} \label{eu_eqn}
{(\Delta +k^2)u^{sca} = 0 \ \text{in}\  D}.
\end{equation}\par
Similarly, the incident field  $u^{inc}$ is also required to satisfy the Helmholtz equation, resulting in the following equation, 
\begin{equation} \label{eu_eqn}
{(\Delta +k^2)u^{inc} = 0\  \text{in}\  D}.
\end{equation}\par
Also,
\begin{equation} \label{eu_eqn}
{u = u^{inc} + u^{sca}}.
\end{equation}
where $u = 0$ on $\partial \Omega$ for a sound soft obstacle, resulting in the following equation,
\begin{equation} \label{eu_eqn}
{(\Delta +k^2)u = 0 \ \text{in} \ D}.
\end{equation}
\subsection{Assumptions and Conditions}
\begin{itemize}
  \item The obstacle is considered to be a sound soft, small, smooth perturbation of a disk.
\item The data used for the reconstruction is of the near field type.
\item The incoming waves can be plane waves, cylindrical waves or point source waves (waves generated by a point source in the exterior of $\Gamma$) and are time harmonic.
\item Only one incident field at a constant frequency and angle is required.
\item $\mathbb{R}^2 \backslash \bar{\Omega}$ is presumed to be filled with a homogenous matter having a fixed wavenumber $k = \frac{2\pi}{\lambda}$.
\item $f(\theta) = \varepsilon g (\theta)$ is assumed.
\item The Helmholtz equation is two dimensional.
\item Only the obstacle surface is reconstructed.
 \end{itemize}
 \section{Reconstruction of shapes and impedance functions using few far-field measurements \cite{RSI}}
 This paper describes a method of reconstruction of the shape and impedance function of complex obstacles \cite{JAM} from few far field acoustic measurements. \par
 The obstacle $D$ is considered to be a star shaped, bounded domain of $\mathbb{R}^2$ such that $\mathbb{R}^2 \backslash \bar{D}$ is connected. The boundary of $D$ is $\partial D$ and its impedance function is given by $\sigma$. Let the media in which propagation is occurring be homogeneous and cylindrical and the fields be time - harmonic and acoustic in nature. \par
 A plane wave $u^i(x,d) = e^{ikdx}$ is incident on the obstacle $D$ and $u^s(x,d)$ be the scattered wave. Let $u^\infty(\hat{x},d)$ be the far field of the scattered wave $u^s(x,d)$ corresponding to the incident direction $d$. Here, $k>0$ is the wave number and $\hat{x} \coloneqq \bigl(\frac{x}{|x|}\bigl)$.\par
 The inverse problem is to find $(\partial D, \sigma)$ given $u^\infty(\hat{x},d)$ for every $\hat{x} \in S^1$ and for $K$ incident directions $d = d_1, d_2, d_3,\ldots,d_K $.
\subsection{Governing equation}
 The governing equation is the Helmholtz equation which is given by, 
 \begin{equation} \label{eu_eqn}
{(\Delta +k^2)u = 0 \ \text{in} \  \mathbb{R}^2\backslash\bar{D}}
\end{equation}
 where $k>0$ is the wave number and $u$ is the total field given by
 \begin{equation} \label{eu_eqn}
{u(x, d)  \coloneqq  u^i(x, d) + u^s(x, d)}.
\end{equation}
 \subsection{Assumptions and Conditions}
 \begin{itemize}
  \item Both the shape and impedance function are reconstructed.
 \item The obstacle is considered to be of a complex structure and is later assumed to be star shaped.
 \item Far field measurements are taken.
 \item The fields are time harmonic and acoustic in nature.
  \item The media in which propagation is occurring is cylindrical and homogeneous. 
 \item It is considered that the boundary $\partial D$ is coated and of class $C^2$.
 \item At the boundary $\partial D$,  it is assumed that $u$ satisfies the Robin type boundary condition, which is given by \ \ \  $\frac{\partial u}{\partial n} + ik\sigma u = 0$ on $\partial D$ where $\sigma$ is the impedance function and $n$ is the outward unit normal of $\partial D$. It is considered that $\sigma$ is a real valued $C^1$-\ continuous function and has a uniform lower bound $\sigma\_>0$ on $\partial D$.
 \item The Sommerfeld radiation condition is satisfied by the scattered field $u^s$,
 \begin{equation} \label{eu_eqn}
{ \lim_{r\to\infty} \sqrt{r} \biggl(\frac{\partial  u^s}{\partial r} - iku^s \biggl) = 0 \ \text{where} \ r = |x|}.
\end{equation}
 \end{itemize}
 \section{Analysis and Conclusion}
 All three methods have provided a unique theoretical approach to solving the inverse scattering problem of the reconstruction of the shape and/or impedance function of an object and have also provided numerical examples that justify their respective methods. Maponi, Recchioni and Zirilli \cite{CMA} and He, Kindermann and Sini \cite{RSI} have reconstructed the shape and impedance of the obstacle whereas Li and Wang \cite{JCP} have reconstructed only the shape.\par
 All three techniques have used the Helmholtz equation as the governing equation. Maponi, Recchioni and Zirilli \cite{CMA} and He, Kindermann and Sini \cite{RSI} consider the shape of the obstacle to be a star and He, Kindermann and Sini \cite{RSI} have an additional condition that there must be knowledge of a point inside the obstacle. Li and Wang \cite{JCP} consider the shape to be a small and smooth perturbation of a ring. The method described by Maponi, Recchioni and Zirilli \cite{CMA} is applicable in three-dimensional geometry whereas the methods described by Li and Wang \cite{JCP} and He, Kindermann and Sini \cite{RSI} are applicable in two-dimensional geometry. Li and Wang \cite{JCP} use near field data for the reconstruction and hence, the reconstruction has subwavelength resolution which is highly desirable. This advantage is not seen in the techniques given by Maponi, Recchioni and Zirilli \cite{CMA} and He, Kindermann and Sini \cite{RSI}, both in which, far field data is used, resulting in lower resolution of the reconstructed data. Li and Wang \cite{JCP} describe a method which is applicable only for sound soft obstacles, a condition which does not arise in the approaches given by Maponi, Recchioni and Zirilli \cite{CMA} and He, Kindermann and Sini \cite{RSI}. In the paper authored by Maponi, Recchioni and Zirilli \cite{CMA}, the method is applicable in the resonance region, that is, when $kL\approx 1$, where $L$ is the characteristic length of the obstacle, a condition which is not required by Li and Wang \cite{JCP} and He, Kindermann and Sini \cite{RSI}. For the numerical experiments, Maponi, Recchioni and Zirilli \cite{CMA} make the use of several incoming waves whereas He, Kindermann and Sini \cite{RSI} make the use of four uniformly distributed incident waves. Li and Wang \cite{JCP} make the use of any type of time harmonic field and only one incident field at a constant frequency and angle is used.\par
 All three methods result in stable, efficient reconstructions. The reconstructed shape is more accurate than that of the impedance function in the techniques described by Maponi, Recchioni and Zirilli \cite{CMA} and He, Kindermann and Sini \cite{RSI}. He, Kindermann and Sini \cite{RSI} give us an approach wherein the reconstruction is of high quality for non-convex obstacles and fixed or dynamic impedance functions. This method works well for initial values of shape and impedance function and substantial noise. \par
 Li and Wang \cite{JCP} list some areas of future study which include the application of the described method to other boundary conditions and problems and the usage of the method for obstacles in three-dimensional geometry. In their paper, the reconstruction of other common shapes and more significant deformations is mentioned for possible further study. Convergence analysis \cite{CAN} of the described algorithm, near field imaging of multiple shapes and multiscale or arbitrary surfaces are also mentioned by them as possible forays in the topic. He, Kindermann and Sini \cite{RSI} mention the local uniqueness question for common shapes and for Robin boundary conditions as a possible direction for further research.\par

\ifCLASSOPTIONcaptionsoff
  \newpage
\fi



%




\bibliographystyle{BibTeXtran}   
\bibliography{BibTeXrefs}       

%








\end{document}